\documentclass[12pt]{amsart}

\usepackage{amsmath, amssymb, graphics}

\newtheorem{theorem}{Theorem}[section]
\newtheorem{lemma}[theorem]{Lemma}
\newtheorem{corollary}[theorem]{Corollary}

\begin{document}

\title[Asymptotic functions]{Asymptotic functions of entire functions}

\author[Aimo Hinkkanen, Joseph Miles, John Rossi]{Aimo Hinkkanen$^{1}$, Joseph Miles$^{2}$, John Rossi$^{3}$}

\address{$^{1}$ $^{2}$ Department of Mathematics, University of Illinois at Urbana--Champaign, 1409 W. Green St., Urbana, IL, 61801, U.S.A.; {aimo@illinois.edu, joe@math.uiuc.edu}}

\address{$^{3}$ Department of Mathematics, Virginia Tech, 225 Stanger St., Blacksburg, VA 24061-1026, U.S.A; rossij@vt.edu}

\dedicatory{Dedicated to the memory of Walter K.~Hayman, FRS}

\date{}


\abstract  
If $f$ is an entire function and $a$ is a complex number, $a$ is said to be an asymptotic value of $f$ if there exists a path $\gamma$ from $0$ to infinity such that $f(z) - a$ tends to $0$ as $z$ tends to infinity along $\gamma$.  The Denjoy--Carleman--Ahlfors Theorem asserts that if $f$ has $n$ distinct asymptotic values, then the rate of growth of $f$ is at least order $n/2$, mean type.  A long-standing problem asks whether this conclusion holds for entire functions having $n$ distinct asymptotic (entire) functions, each of growth at most order $1/2$, minimal type. In this paper conditions on the function $f$ and associated asymptotic paths are obtained that are sufficient to guarantee that $f$ satisfies the conclusion of the Denjoy--Carleman--Ahlfors Theorem.  In addition, for each positive integer $n$, an example is given of an entire function of order $n$ having $n$ distinct, prescribed asymptotic functions, each of order less than $1/2$.
\endabstract

\maketitle

\section{Introduction}

Suppose that $f(z)$ is an entire function in the complex plane ${\mathbb C}$. An entire function $a(z)$ is said to be an asymptotic function for $f$ if there exists a path $\gamma$ in ${\mathbb C}$ from $0$ to infinity such that $f(z)-a(z)$ tends to $0$ as $z$ tends to infinity along $\gamma$; we then say that $f$ is asymptotic to $a$ on $\gamma$. In the case that $a(z)$ is a constant function, it is said to be an asymptotic value of $f$. The Denjoy--Carleman--Ahlfors Theorem (\cite{A1}; in a very general form for subharmonic functions this can be found in \cite{Hsub2}, Theorem~8.9, p.~562) states that if $f$ has $n$ distinct asymptotic values then the rate of growth of $f$ is at least order $n/2$, mean type. This bound is known to be sharp (\cite{A}, p.~1). For consider $f(z) = \int_0^z (w^{-n/2} \sin w^{n/2})\, dw$. Then $f$ has order $n/2$, and if $0\leq \nu \leq n-1$ and $z=r e^{ 2\nu i \pi/n }$, where $r>0$, we have
$$f(z)\to e^{ 2\nu i \pi/n } \int_0^{\infty}  x^{-n/2} \sin ( x^{n/2} )  \, dx$$ as $r=|z|\to \infty$. 

It is an open question \cite{Hprob}, Problem~2.3, whether the analogue of the Denjoy--Carleman--Ahlfors Theorem  holds for asymptotic functions $a(z)$ with growth at most order $1/2$, minimal type. For such functions it is known that the minimum modulus is unbounded on a sequence of circles with radii tending to infinity. Absent this minimum modulus condition on the asymptotic functions $a(z)$, there are easy examples of entire functions of finite order with infinitely many asymptotic functions. Such an example is $f(z)=\frac{ \sin \sqrt{z} }  {   \sqrt{z} } $ or $f(z)=e^{-z}$, and $a_j(z)=jf(z)$ for each positive integer $j$. Note that both $f(z)$ and $a_j(z)$ tend to $0$ as $z\to\infty$ along the positive real axis.

Throughout this paper we let $B(0,r)$ denote the disk $\{z\in {\mathbb C} \colon |z|<r \}$ and let $S(0,r)$ be its boundary circle $\{z\in {\mathbb C} \colon |z|=r \}$, where $r>0$. If $f$ is entire and $D$ is a domain in ${\mathbb C}$, we write
$$
M(r,D,f) = \sup \{ |f(z)| \colon z\in D\cap S(0,r) \} .
$$
We denote the usual maximum modulus $M(r,{\mathbb C},f)$ of $f$ by $M(r,f)$. A path $\gamma$ from $0$ to infinity is said to be segmental if it is a Jordan arc consisting of line segments whose vertices accumulate only at infinity.

Recall that the order $\rho=\rho(f)$ of an entire function $f$ is defined by 
$$
\rho(f) = \limsup_{r\to\infty} \frac{ \log \log M(r,f) } { \log r } .
$$
If $0<\rho<+\infty$, we say that $f$ is of minimal type if
$\log M(r,f) = o(r^{\rho})$ as $r\to\infty$, and we say that $f$ is of mean type if 
$$
0<   \limsup_{r\to\infty} \frac{  \log M(r,f) } { r^{\rho} }      < +\infty .
$$

The strongest result on Problem 2.3 in \cite{Hprob} is due to Fenton \cite{F}, who proved the following theorem.

{\bf Theorem A.} {\sl If $f$ is entire with $n$ distinct asymptotic functions of order less than $1/4$, then the growth of $f$ is at least order $n/2$, mean type.}

In \cite{Hsub2}, Theorem~8.13, p.~577, Hayman showed that Fenton's methods can be used to obtain this result for asymptotic functions with growth no faster than order $1/4$, minimal type. For an earlier result of Somorjai with asymptotic functions of order $<1/30$, see \cite{S}. 

Stronger results are known if the asymptotic paths are rays. In this case Denjoy (\cite{D1}, \cite{D2}) showed that if $f$ has order $\mu$ then the number of distinct asymptotic functions with order less than $1/( 2 + \mu^{-1}  )$ is at most $2\mu$.

Dudley Ward and Fenton \cite{DWF} obtained the following result.

{\bf Theorem B.} {\sl Suppose that $f$ is analytic in the sector $D=\{ z\in {\mathbb C} \colon |\arg z|<\eta \}$ for some $\eta\in (0,\pi)$ and is continuous on $\partial D$.    
Suppose that $a(z)$ and $b(z)$ are entire, each with order less than 
$1/(2 + 2\eta \pi^{-1})$,  are not identically zero and satisfy
\begin{equation} \label{11}
f( t e^{i\eta} ) - a ( t e^{i\eta} )  \to 0\quad \text{ as }\, t\to + \infty 
\end{equation}
and
\begin{equation} \label{12}
f( t e^{ - i \eta} ) - b ( t e^{ - i \eta} )  \to 0\quad \text{ as }\, t\to + \infty . 
\end{equation}
Then
$$
\liminf_{r\to\infty} \frac{   \log M(r,D,f)   } {   r^{   \pi/(2\eta)    }     } > 0  .
$$
}

Elementary arguments lead from Theorem B to the following corollary.

{\bf Corollary A.} {\sl Suppose that $f$ is entire and that $\gamma_1,\gamma_2, \dots , \gamma_n$ are distinct rays emanating from $0$. Suppose for $1\leq j\leq n$ that $a_j(z)$ are distinct entire functions of order less than $1/(2 + 2 n^{-1} ) $ and that $f(z)$ is asymptotic to $a_j(z)$ on $\gamma_j$. Then
$$
\liminf_{r\to\infty} \frac{   \log M(r,f)   } {   r^{   n/2    }     } > 0  .
$$
}

Hinkkanen and Rossi \cite{HR} obtained the conclusion of Theorem B under the relaxed assumption that the open set $D$ is bounded by paths $\gamma_1$ and $\gamma_2$, not necessarily disjoint away from $0$, with (\ref{11}) holding on $\gamma_1$ and (\ref{12}) holding on $\gamma_2$ and that the angular measure $m ( D \cap S(0,r) )$ of $D \cap S(0,r)$ satisfies 
$$
m ( D \cap S(0,r) ) \leq 2 \eta 
$$
for all large $r$. 

For further results in this direction, see \cite{FR1}, \cite{FRRD1}, \cite{FRRD2},   \cite{FRRD3}, and \cite{FR2}.

The principal result of this paper is the following theorem.

\begin{theorem} \label{th1}
Suppose that $\gamma_1,\gamma_2, \dots , \gamma_n$ are $n$ simple segmental paths from $0$ to infinity, disjoint except at the origin, arranged in counterclockwise fashion. Let $D_j$ be the Jordan domain between $\gamma_j$ and $\gamma_{j+1}$, where $\gamma_{n+1}=\gamma_1$. Suppose that $f$ is entire and that $a_j(z)$, for $1\leq j\leq n$, are distinct entire functions with order $\rho(a_j)<1/2$ such that $f$ is asymptotic to $a_j$ on $\gamma_j$. Suppose that there exists a number
$\kappa > \rho := \max \{ \rho(a_j) \colon 1\leq j\leq n \}$ such that
\begin{equation} \label{13}
\limsup_{z\to\infty \atop z\in D_j}  \frac{  \log |f(z)|  } {   |z|^{\kappa}     } > 0
\end{equation}
for each $j$. Then 
\begin{equation} \label{14}
\liminf_{r\to\infty}  \frac{  \log M(r,f)   } {   r^{n/2}      } > 0  .
\end{equation}
\end{theorem} 

We note in particular that if (\ref{13}) holds with $\kappa=1/2$, then (\ref{14}) holds.

 From Theorem~\ref{th1} we obtain the following corollary.
 
 \begin{corollary} \label{co1}
Suppose for $1\leq j\leq n$ that $a_j(z)$ are distinct entire functions with order $\rho(a_j)<1/2$. For $1\leq j\leq n$, let $\gamma_j$ be the ray with $\arg z = 2\pi j/n$. Suppose that $f$ is entire and that $f$ is asymptotic to $a_j$ on $\gamma_j$. Then 
$$
\liminf_{r\to\infty}  \frac{  \log M(r,f)   } {   r^{n/2}      } > 0  .
$$
 \end{corollary}
 
In comparing Corollary~\ref{co1} with Corollary A, we note that Corollary~\ref{co1} treats asymptotic functions $a(z)$ of all orders less than $1/2$, but requires that the rays be equally spaced. 

We also show  
that Fenton's result in the special case that the $a_j(z)$ are polynomials follows quite easily from Theorem~\ref{th1}. 

Finally, in Section~\ref{ex} we give an example of an entire function of order $n$ having $n$ distinct, prescribed asymptotic functions of order $<1/2$. Our method does not seem capable of constructing an entire $f$ of order $n$ with $2n$ such prescribed  asymptotic functions. 

\section{Proofs}

\subsection{A lemma}

Fundamental to our approach is the following lemma.

\begin{lemma} \label{le1}
Suppose that $\gamma_1$ and $\gamma_2$ are simple segmental paths from $0$ to infinity, disjoint except for the origin. Let $D$ be a Jordan domain with $\partial D = \gamma_1 \cup \gamma_2$. For $t>0$, let $\Phi(t)$ be the angular measure of $D\cap S(0,t)$. Suppose that $a_1$ and $a_2$ are distinct entire functions with orders $\rho(a_1)$ and $\rho(a_2)$ satisfying
$\rho = \max \{    \rho(a_1), \rho(a_2)       \} < 1/2$, 
such that $f$ is asymptotic to $a_j$ on $\gamma_j$ for $j=1,2$. 
Suppose that there exists a number $\kappa>\rho$ such that 
$$
\limsup_{z\to\infty \atop z\in D}  \frac{  \log |f(z)|  } {   |z|^{\kappa}     } > 0  . 
$$
Then there exists $R_1>0$ such that for all $R>R_1$, we have
$$
\log M(R,D,f) \geq \frac{ \pi } { 8 } \exp \left\{ \pi \int_{R_1}^R \frac{  dt } {  t \Phi(t) }        \right\}
\geq   \frac{ \pi } { 8 }  \left(   \frac{ R } { R_1 }     \right)^{1/2}  . 
$$
\end{lemma}

{\bf Proof.} Without loss of generality, we may assume that $\kappa<1/2$. Select numbers $\kappa_1$ and $\kappa_2$ such that $\rho< \kappa_1 < \kappa_2 < \kappa$. There exists a number $C\geq 1$ such that for all $z\in \partial D$ we have
$$
\max \{  \log |a_1(z)| ,   \log |a_2(z)|      \} < C + |z|^{\kappa_1}  .
$$
There exists $M' \geq 1$ such that for $j=1,2$, 
$$
\log | f(z) - a_j(z) | \leq M' 
$$
for all $z\in \gamma_j$ and hence
$$
\log | f(z) |  \leq \log^+ | a_j(z) | + M' + \log 2 < C + |z|^{\kappa_1} + M' + \log 2 .
$$

With $M = C + M' + \log 2$, we have
\begin{equation} \label{21}
\log | f(z) | < M +  |z|^{\kappa_1} 
\end{equation}
for all $z\in \partial D$. 

Set
$$
A_0 = \frac{  20    } {  \left(  \frac{1}{2} - \kappa_1     \right) 4^{   \frac{1}{2} - \kappa_1     }        }  > 1 .
$$

We choose $R_0$ so large that for all $r\geq R_0$ we have
$$
1 + M + A_0 r^{\kappa_1} < r^{\kappa_2}  .
$$

The following lemma is from \cite{Hsub2}, Lemma~8.13, p.~583.  

\begin{lemma} \label{hay1}
Let $\phi$ be a non-negative continuous  convex function of $\log t$ for $0\leq t<\infty$ with $\phi(0)=0$, and suppose that for some $\delta>0$, we have
$$
\int_{\delta}^{\infty} \frac{ \phi(t) } { t^{3/2}  } \, dt < \infty .
$$
Let $D$ be a  domain in the plane such that every boundary point of $D$ is regular for Dirichlet's problem and such that for all $r\in (0,\infty)$, the circle $S(0,r)$ intersects the complement of $D$. 
Then there is a function $u$, continuous and non-negative in $\overline{D}$ and harmonic in $D$, such that $u(z)=\phi(|z|)$ for all $z\in \partial D$ and such that for all $z\in D$, we have
\begin{equation} \label{a11}
\phi(|z|) \leq u(z) \leq 20 |z|^{1/2}  \int_{4|z|}^{\infty}  \frac{ \phi(t) } { t^{3/2}  } \, dt  .
\end{equation}
\end{lemma}

Hayman used an estimate for harmonic measure that is not as strong as is known. Hence the number $20$ on the right hand side of (\ref{a11}) can be reduced but we do not try to do this as it is not important for us.

We use $\phi(t) = t^{\kappa_1}$ for our earlier choice of $\kappa_1$ and calculate
$$
20 |z|^{1/2}  \int_{4|z|}^{\infty}  \frac{ \phi(t) } { t^{3/2}  } \, dt 
=  A_0   |z|^{\kappa_1}   .
$$
We apply Lemma~\ref{hay1} 
to obtain a harmonic function $u$ on $D$ such that for all $z\in D$ we have
$$
 |z|^{\kappa_1}  \leq u(z) \leq A_0  |z|^{\kappa_1}   .
$$

By our assumptions, there exists $z_1\in D$ with $|z_1|>R_0$ such that
\begin{equation} \label{22}
\log | f(z_1) | > |z_1|^{\kappa_2} > A_0 |z_1|^{\kappa_1} + M + 1 .
\end{equation}
We set $R_1=|z_1|$. 

Suppose that $R>R_1$. Let $D(R)$ be the component of $D\cap B(0,R)$ containing $z_1$. We note that $\partial D(R)$ lies in $\partial D$ except for a non-empty union of open arcs in $S(0,R)$. Let $U=U_R$ be the harmonic function in $D(R)$ that satisfies $U(z)=0$ if $z\in \partial D(R)$ lies in the interior of an arc in $S(0,R)$ lying in $\partial D(R)$, while $U(z)=|z|^{\kappa_1}$ if $z\in \partial D(R)$ and $|z|<R$ (and thus $z\in \partial D$). Let $u$ be the harmonic function obtained above from Lemma~\ref{hay1} with the choice $\phi(t)=t^{\kappa_1}$. Evidently $U-u$ is harmonic on $D(R)$ and $U(z)-u(z)\leq 0$ for all $z\in \partial D(R)$. Thus $U(z)-u(z) \leq 0$ for all $z\in D(R)$. 

Let
$$
M'(R,D,f) = \max \{ |f(z) | \colon |z|=R, \,\, z\in \partial D(R) \} \leq M(R,D,f) .
$$ 
Write $H(R) = S(0,R) \cap \partial D(R)$. Let $\omega(R,z)$ denote the harmonic measure of $H(R)$ at $z\in D(R)$. The function
$$
w(z) =  U(z) + M + \omega(R,z) \log M'(R,D,f) 
$$
is harmonic on $D(R)$ and
$$
\log |f(z)| \leq w(z) 
$$
for all $z\in \partial D(R) $ by (\ref{21}). Then for all $z\in D(R)$ 
\begin{equation} \label{23}
\log | f(z) | \leq w(z) \leq u(z) + M + \omega(R,z) \log M(R,D,f) 
.
\end{equation}
Setting $z=z_1$ we obtain from (\ref{22})
$$
A_0 |z_1|^{\kappa_1} + M + 1 \leq \log | f(z_1) | \leq 
A_0 |z_1|^{\kappa_1} + M + \omega(R,z_1) \log M(R,D,f)  , 
$$
implying that
$$
1 \leq \omega(R,z_1) \log M(R,D,f)  .
$$

If we now let $\Phi^* (t)$ be the angular measure of $D(R) \cap S(0,t)$, we have
$\Phi^* (t) \leq  \Phi (t)$, and by \cite{GM}, Theorem~6.2, (6.4), p.~149 and p.~158, we get 
\begin{equation} \label{24}
\omega(R,z_1) \leq 
\frac{ 8 } { \pi  } \exp \left\{ - \pi  \int_{R_1}^R \frac{  dt } {  t \Phi^*(t) }        \right\} \leq 
\frac{ 8 } { \pi  } \exp \left\{ - \pi  \int_{R_1}^R \frac{  dt } {  t \Phi(t) }        \right\}
.
\end{equation}
We rearrange to obtain
$$
\log M(R,D,f) \geq \frac { \pi  } { 8 } \exp \left\{  \pi  \int_{R_1}^R \frac{  dt } {  t \Phi(t) }        \right\} .
$$
The second inequality in our conclusion follows from the fact that $\Phi(t) \leq 2\pi$. 

\subsection{Proof of Theorem~\ref{th1}}

We now turn to the proof of Theorem~\ref{th1}. For each $j$, let $\Phi_j(t)$ be the angular measure of $D_j\cap S(0,t)$. We apply Lemma~\ref{le1} on each $D_j$ to conclude that there exists $R_1(j)>0$ such that for all $R>R_1(j)$, 
$$
\log M(R,D_j,f) \geq \frac { \pi  } { 8 } \exp \left\{  \pi  \int_{R_1(j)}^R \frac{  dt } {  t \Phi_j(t) }        \right\} .
$$
Let $R_1 = \max\{ R_1(j) \colon 1\leq j\leq n \}$. Then for all $R>R_1$ we have
$$
\log M(R,f) \geq \frac { \pi  } { 8 } \exp \left\{  \pi  \int_{R_1(j)}^R \frac{  dt } {  t \Phi_j(t) }        \right\}
$$
for each $j$ with $1\leq j\leq n$. 

We have
$$
n^2 \leq 
\left( \sum_{j=1}^n  \frac{1}{ \Phi_j(t) }    \right)  \left( \sum_{j=1}^n   \Phi_j(t)   \right)  
\leq 
 2\pi \sum_{j=1}^n  \frac{1}{ \Phi_j(t) }  ,
$$
implying that
$$
\frac{ n^2 } { 2 } \int_{R_1}^R \frac{ dt } { t } \leq \pi \sum_{j=1}^n   \int_{R_1}^R 
\frac{ dt } { t  \Phi_j(t) } .
$$
Thus for all $R>R_1$, there exists $j=j(R)$ such that
$$
\frac{ n } { 2 } \int_{R_1}^R \frac{ dt } { t } \leq \pi    \int_{R_1}^R 
\frac{ dt } { t  \Phi_j(t) } .
$$
So for $R>R_1$, with $j=j(R)$ we have
$$
\log M(R,f) \geq \frac { \pi  } { 8 } \exp \left\{  \pi  \int_{R_1(j)}^R \frac{  dt } {  t \Phi_j(t) }        \right\}
\geq \frac { \pi  } { 8 } \left(  \frac{ R } { R_1 }     \right)^{n/2} ,
$$
establishing ({\ref{14}). 

\subsection{Proof of Corollary~\ref{co1}}

We now prove Corollary~\ref{co1}. Let $D_j$ be the sector between $\gamma_j$ and $\gamma_{j+1}$ (with $\gamma_{n+1}=\gamma_1$). From elementary considerations we have
$$
\limsup_{r\to\infty} \frac{ \log M(r,f) } { r^{1/2} } > 0 .
$$
Trivially for each $R$ there exists $j=j(R)$ such that 
$ \log M(R,f) = \log M(R,D_j,f)$. Thus there exists $j$ such that
$$
\limsup_{z\to\infty \atop z\in D_j} \frac{ \log |f(z)| } { |z|^{1/2} } > 0 .
$$
We may thus apply Lemma~\ref{le1} on this $D_j$ to conclude that there exists $R_1=R_1(j)>0$ such that for all $R>R_1$
$$
 \log M(R,f) \geq \log M(R,D_j,f) \geq \frac { \pi  } { 8 }  
\exp \left\{  \pi  \int_{R_1(j)}^R \frac{  dt } { 2\pi t n^{-1} }        \right\}
=  \frac { \pi  } { 8 }    \left(  \frac{ R } { R_1 }     \right)^{n/2}  .
$$
\qed

{\bf Remark.}  Our proofs in fact yield somewhat stronger versions of Theorem~\ref{th1} and Corollary~\ref{co1}.  The assumption that all of the asymptotic functions are distinct can be weakened.  Both Theorem~\ref{th1} and Corollary~\ref{co1} depend on applications of Lemma~\ref{le1}, each application requiring only that the asymptotic functions associated with adjacent paths are distinct. Thus both results hold even if $f$ has fewer than $n$ distinct asymptotic functions on the $n$ disjoint paths, provided only that asymptotic functions associated with adjacent paths are distinct.

Similar considerations apply to the Denjoy--Carleman--Ahlfors Theorem and Fenton's theorem, as an examination of the arguments in \cite{A1} and \cite{F} shows.  In fact, more can be said in the case of asymptotic values.  A path to infinity on which $f$ tends to a constant $a$, finite or infinite, corresponds to a singularity of $f^{-1}$ lying over $a$.  More than one singularity can lie over a given $a$.  Certain singularities are classified  as direct singularities.  If $f$ is meromorphic of order $\rho$, Ahlfors \cite{A2} showed that $f^{-1}$ has no more than $\max \{ 2 \rho , 1\}$ direct singularities.  If $f$ is entire, all singularities lying over infinity are direct.  From this it follows that if $f$ is entire of order $\rho$, $f$ has no more than $2\rho$ singularities lying over finite points, implying the Denjoy--Carleman--Ahlfors Theorem. For details, including a necessary modification of Ahlfors's arguments, see   \cite{N36}, pp.~309--313 or  \cite{N}, pp.~303--307.

\subsection{Fenton's theorem for polynomials $a_j$}

We next turn to a proof, based on Theorem~\ref{th1}, of Fenton's theorem in the case that all $a_j$ are polynomials.  We treat the situation where not all $a_j$ are assumed to be distinct.

Suppose that $f$ is entire and that ${\mathcal F} = \{ \gamma_1, \gamma_2,  \dots ,\gamma_m\}$ is a collection of paths from $0$ to infinity, each with an associated polynomial $a_j$ such that $f$ is asymptotic to $a_j$ on $\gamma_j$.  Suppose that  two of these paths, say $\gamma_j$ and $\gamma_k$, intersect on a sequence of points tending to infinity.  Then $a_j - a_k$ tends to $0$ on this sequence, implying that $a_j\equiv a_k$;  in this circumstance, we delete one of $\gamma_j$ and $\gamma_k$, obtaining a subcollection of ${\mathcal F}$.  If two paths in this subcollection intersect on a sequence tending to infinity, we again delete one of them.  We continue until we arrive at a subcollection ${\mathcal F}'$ with the property that there is a disk $D$ centered at the origin containing all points of intersection of any two members of ${\mathcal F}'$.  After renumbering, we write ${\mathcal F}' = 
\{\gamma_1, \gamma_2,  \dots , \gamma_q\}$.  After an obvious modification of the paths on $D$, we may suppose that the paths in ${\mathcal F}'$ are simple, segmental, disjoint except at the origin, and numbered in counterclockwise order.  If any two adjacent paths in ${\mathcal F}'$, say $\gamma_j$ and $\gamma_{j+1}$ (with $\gamma_{q+1} = \gamma_1$), are associated with the same asymptotic polynomial, we delete one of these paths.  We continue this process until we obtain a subcollection ${\mathcal F}'' = \{\gamma_1, \gamma_2,  \dots , \gamma _n\}$ of disjoint simple segmental paths numbered in counterclockwise order such that the polynomials associated with any two adjacent paths are distinct. 

Every collection ${\mathcal F}$ has such a subcollection ${\mathcal F}''$.  Our task is to show that $f$ has growth at least order $n/2$, mean type, where $n$ is the number of paths in such a subcollection ${\mathcal F}''$.  We note that if all the asymptotic polynomials associated with the paths in ${\mathcal F}$ are distinct, then ${\mathcal F} = {\mathcal F}''$, after obvious modifications of the paths on a disk $D$.
Similarly, ${\mathcal F} = {\mathcal F}''$ (after obvious modifications) if all the paths in ${\mathcal F}$ are disjoint outside some disk and the asymptotic polynomials associated with adjacent paths are distinct.

We adopt the notation of Theorem~\ref{th1}.  We consider such a collection ${\mathcal F}''$ and for $1 \leq j \leq n$ write
$$
a_j(z) = \sum_{p=0}^{d_j} b_{pj} z^p  
$$
where $d_j$ is the degree of $a_j$. 

Fix $ j$ and assume, without loss of generality, that $d_j \geq d_{j+1}$.  Replacing $f$ by $ f - a_{j+1}$ and $a_j$ by $a_j - a_{j+1}$ without changing notation, we may assume $ a_{j+1} \equiv 0$ and $a_j \not\equiv 0$.  Our goal is to show that in $D_j$, the modified function $\log |f(z)|$ grows at least like $|z|^{1/2}$, clearly implying that our original $\log |f(z)|$ does as well.

If $a_j$ is constant, it follows from familiar arguments in the case of asymptotic values that  $f$ is unbounded in $D_j$.  On $\partial D_j$, we have $\log |f(z)| \leq M$ for some $M>1$. For every $z\in \partial D_j$ with $|z|<R$, we have
(\ref{23}) with $u(z) =0$. We note that
$$
\omega(R,z) = O\left(  \frac{ |z| } { R }     \right)^{1/2}  
$$
from (\ref{24}). 
Letting $M_j(R)=M(R,D_j,f)$, we find from (\ref{23}) that if there is a sequence of $R\to\infty$ on which $\log M_j(R)=o(R^{1/2})$, then $f$ is bounded in $D_j$, which is a contradiction. Thus we must have 

\begin{equation} \label{a35} 
\liminf_{R\to\infty} \frac{ \log M_j(R)   } {  R^{1/2}  } >0 .
\end{equation}

Suppose now that $d_j\geq 1$. Let $\alpha$ be a value that $f$ takes infinitely often, say (at least) at distinct values $c_p$  for $1\leq p\leq d_j$. Set $P(z) = \prod_{p=1}^{d_j} (z-c_p)$ and
$$
g(z) = \frac{ f(z) - \alpha  } {  P(z) }  .
$$
Then $g$ is entire, and $g(z)\to 0$ as $z\to\infty$ along $\gamma_{j+1}$ (since $f(z)\to 0$). As $z\to\infty$ along $\gamma_{j}$, we have $f(z) - a_j(z) \to 0$, so that
$$
g(z) = o(1) +  \frac{ a_j(z) - \alpha  } {  P(z) } = o(1) + \frac{ a_j(z)   } {  P(z) }
= o(1) +   b_{d_j j}   
$$
since $\deg a_j = \deg P = d_j$. Recall that $b_{d_j j}\not= 0$. Applying the previous arguments to $g$ instead of $f$, we see that (\ref{a35}) holds if $M_j(R)$ refers to $g$, and hence also holds if $M_j(R)$ refers to $f$. The above works in all domains $D_j$, so that we may apply Theorem~\ref{th1} to conclude that 
$$
\liminf_{R\to\infty} \frac{ \log M(R,f)   } {  R^{n/2}  } >0 .
$$

\section{An example} \label{ex}

\begin{theorem} \label{th2}
Let $n$ be a positive integer. For $1\leq j\leq n$, suppose that $a_j$ is an entire function of growth no faster than order $n$, minimal type. There exists an entire function $f$ of order $n$, mean type, such that each $a_j$, for $1\leq j\leq n$, is an asymptotic function of $f$.
\end{theorem}

Our proof is an adaptation of a technique that Fuchs and Hayman (\cite{FH}, see also \cite{Hmero}, pp.~80--83) used to assign deficiencies of an entire function arbitrarily subject only to the condition that the sum of the deficiencies  does not exceed $1$. We prove the following lemma.

\begin{lemma} \label{le-2}
Let $n$ be a positive integer. Let
$$
c_n = \int_0^{\infty} e^{-t^n} \, dt .
$$
Let the path $\Gamma$ be the boundary of a sector of opening $\frac { 2\pi } { n } $ given by
$$
\Gamma(t) = \left\{  \begin{matrix}  -t e^{-i\pi /n}, & -\infty < t < 0 , \\
t e^{ i\pi /n}, & 0 \leq t < \infty  .    \end{matrix}   \right.
$$
Let 
$\Omega_1 = \{ r e^{i\theta} \colon r>0 \,\, \text{ and } \,\, |\theta|< \pi/n \} $
be the inside of $\Gamma$ and 
$\Omega_2 = \{ r e^{i\theta} \colon r>0 \,\, \text{ and } \,\, 
\pi/n < \theta < 2\pi - \pi/n
 \} $
be the outside. Then there exists an entire function $\varphi$ such that
\begin{equation} \label{31}
\varphi(z) = \left\{  \begin{matrix}   e^{z^n} - \left( \frac{c_n}{\pi} \sin \frac{\pi}{n} \right) z^{-1} + O( |z|^{-2} )    , & z\in \Omega_1 , \\
- \left(  \frac{c_n}{\pi} \sin \frac{\pi}{n} \right)  z^{-1} + O( |z|^{-2} ) , & z\in \Omega_2  .    \end{matrix}   \right.
\end{equation}
uniformly as $z\to\infty$.
\end{lemma}

{\bf Proof of Lemma~\ref{le-2}.} For $z\in {\mathbb C}\setminus \Gamma$, define
$$
E(z) = \frac{ 1 } { 2\pi i} \int_{\Gamma} \frac{ e^{w^n}   } {  w-z    } \, dw  .
$$
We note that if $w\in \Gamma$, then $w^n=-|w|^n$ and thus $E$ is analytic for $z\in  {\mathbb C} \setminus \Gamma$. 
Let
$$
I =  \frac{ 1 } { 2\pi i} \int_{\Gamma}   e^{w^n}    \, dw = 
\frac{ 1 } { 2\pi i} \int_0^{\infty}  e^{ - t^n} \left(  e^{ i\pi /n} - e^{-i\pi /n}     \right) \, dt
= \frac{ 1 } { \pi } c_n \sin \frac{\pi}{n} . 
$$
Note that
$$
\frac{ 1 } { w-z } + \frac{ 1 } { z } = \frac{ w } { z ( w - z ) }  .
$$
Thus
\begin{equation} \label{33}
L := \frac{ 1 } { 2\pi i} \int_{\Gamma} \frac{ w e^{w^n}   } {  z ( w-z )   } \, dw
= E(z) + \left(  \frac{ c_n } { \pi }  \sin \frac{\pi}{n}    \right) \frac{ 1 } { z }    .
\end{equation}

Consider $z\notin \Gamma$ with $|z|=R$ for large $R > 0$. Let $\Gamma_1$ be the portion of $\Gamma$ with modulus at most $R/2$, $\Gamma_2$ the portion of $\Gamma$ with modulus between $R/2$ and $2R$, and $\Gamma_3$ the portion of $\Gamma$ with modulus at least $2R$. Let
$$
d_n = \int_0^{\infty} t e^{ - t^n } \, dt  .
$$
We have
\begin{equation} \label{34}
\left|   \frac{ 1 } { 2\pi i} \int_{\Gamma_1}  \frac{ w } { z ( w - z ) }\,   e^{w^n}    \, dw  \right| + 
\left|   \frac{ 1 } { 2\pi i} \int_{\Gamma_3}  \frac{ w } { z ( w - z ) }\,  e^{w^n}    \, dw   \right| 
< \left(  \frac{   2 } { \pi R^2    } +    \frac{   1 } { \pi R^2    }   \right) d_n   .
\end{equation}

If the distance from $z$ to $\Gamma_2$ is at least $1$, we set $\Gamma_2^* = \Gamma_2$. If the distance from $z$ to $\Gamma_2$ is less than $1$, we indent away from $z$ the segment  of $\Gamma_2$ with distance less than $1$ from $z$ to the arc of circle with distance at least $1$ from $z$ to obtain a modified  $\Gamma_2^* $. 

By Cauchy's Theorem we have
$$
\int_{\Gamma_2}  \frac{ w e^{w^n}  } { z ( w - z ) }     \, dw
= \int_{\Gamma_2^*}  \frac{ w e^{w^n}  } { z ( w - z ) }     \, dw  .
$$
We conclude that
\begin{equation} \label{35}
\left|   \frac{ 1 } { 2\pi i} \int_{\Gamma_2}   \frac{ w e^{w^n}  } { z ( w - z ) }   \, dw  \right|  
= O\left(  R e^{ - (R/2)^n   }     \right)
= O\left(    \frac{ 1 } { R^2 }   \right)
   .
\end{equation}

Combining (\ref{33}), (\ref{34}), and (\ref{35}), we obtain
\begin{equation} \label{36}
\left|  E(z) + \left( \frac{c_n}{\pi} \sin \frac{ \pi } {n }      \right)  \frac{ 1 } { z }  \right|  
= O\left(    \frac{ 1 } { R^2 }   \right)
\end{equation}
for $z\notin \Gamma$. 

Let $E_2(z)$ be $E(z)$ for $z\in \Omega_2$ and let $E_1(z)$ be $E(z)$ in $\Omega_1$. By replacing a short segment of $\Gamma$ by a circular arc lying in $\Omega_1$, we see from Cauchy's Theorem that $E_2$ can be continued analytically across $\Gamma$ onto $\Omega_1$ and, by the Cauchy Integral Formula, that $E_2$ is entire and
\begin{equation} \label{37}
E_2(z) = E_1(z) + e^{z^n} 
\end{equation}
for all $z\in \Omega_1$. We set $\varphi(z) = E_2(z)$ and conclude from (\ref{36}) and (\ref{37})
that  (\ref{31})  holds. This proves Lemma~\ref{le-2}.

We now turn to the proof of Theorem~\ref{th2}. 
Letting $\varphi$ be as in Lemma~\ref{le-2} for $1\leq j\leq n$ we set
$$
\varphi_j(z) = \varphi \left(  e^{   -2\pi i  j / n   }    z  \right)   .
$$
Define
\begin{equation} \label{38}
f(z) = \sum_{j=1}^n \frac{  \varphi_j(z)  a_j(z)    } {   e^{ z^n }        }  .
\end{equation}
Clearly $f$ is entire of order $n$ mean type.

Let $\gamma_j$ be the ray $\{ r e^{  2\pi ij/n      }    \colon r> 0 \}$. Suppose that $z\in \gamma_{j_0}$. Then $  e^{   -2\pi i  j_0 /n    }    z     >0$ and so lies in $\Omega_1$. Thus
$$ \varphi_{j_0}(z) = \varphi\left(  e^{   -2\pi i  j_0 /n   }    z     \right) = 
e^{z^n} + O( |z|^{-1} ) = e^{ |z|^n } + O( |z|^{-1} )  .
$$
Hence
$$
\frac{  \varphi_{j_0}(z)  a_{j_0}(z)    } {   e^{ z^n }        } - a_{j_0}(z) 
= \left(  \frac{ e^{z^n}  + O( |z|^{-1} )  } {  e^{z^n}    }       -1 \right) a_{j_0}(z) 
= \frac{  O( |z|^{-1} )  a_{j_0}(z)  }  {  e^{z^n}     } 
,
$$
implying that
\begin{equation} \label{39}
\lim_{z\to\infty \atop z\in \gamma_{j_0} }
\left(   \frac{  \varphi_{j_0}(z)  a_{j_0}(z)  }  {  e^{z^n}     } -  a_{j_0}(z)  \right)   =0
 .
\end{equation}
We now consider $z\in \gamma_{j_0}$ and $j\not= j_0$. We note that
$ e^{  - 2\pi ij /n  } z     \in \Omega_2$. By (\ref{31}), 
$$
| \varphi_j(z) | = \left|      \varphi \left(  e^{   -2\pi i  j /n   }    z  \right)         \right| 
=  O( |z|^{-1} )  .
$$
Consequently, for $z\in \gamma_{j_0}$,
\begin{equation} \label{310}
\left|  \displaystyle\sum_{j=1 \atop j\not= j_0}^n \frac{  \varphi_j(z)  a_j(z)    } {   e^{ z^n }        }      \right| \leq 
\frac{ \left(   \displaystyle\sum_{j=1 \atop j\not= j_0}^n | a_j(z) |      \right)   O( |z|^{-1} )  } {   e^{ |z|^n }     } = o(1) .
\end{equation}

The combination of (\ref{38}), (\ref{39}), and (\ref{310}) 
completes the proof of Theorem~\ref{th2}.

\end{document}